\documentclass{amsart}

 \theoremstyle{definition}
 
 \theoremstyle{remark}
 
 \numberwithin{equation}{subsection}

\begin{document}

\centerline{\bf 2-Local derivations on algebras of matrix-valued
functions on a compact}

\medskip

\bigskip

\centerline{Ayupov Shavkat Abdullayevich}

\medskip

\centerline{Institute of Mathematics, National University of Uzbekistan, Tashkent, Uzbekistan}

\centerline{sh$_-$ayupov@mail.ru}

\bigskip

\centerline{Arzikulov Farhodjon Nematjonovich (a corresponding
author)}

\medskip

\centerline{Faculty of Mathematics, Andizhan State University,
Andizhan, Uzbekistan}

\centerline{arzikulovfn@rambler.ru}

\medskip

\centerline{\bf Abstract}

In the present paper 2-local derivations on various algebras of
infinite dimensional matrix-valued functions on a compact are
considered. It is proved that every 2-local derivation on such
algebra is a derivation. Also we explain that the method developed
in the given paper can be applied to associative, Jordan and Lie
algebras of infinite dimensional matrix-valued functions on a
compact.

\medskip


\medskip

{\bf AMS 2000:} 46L57, 46L40

\medskip

\section*{Introduction}

\medskip

The present paper is devoted to 2-local derivations on algebras.
Recall that a 2-local derivation is defined as follows: given an
algebra $A$, a map $\bigtriangleup : A \to A$ (not linear in
general) is called a 2-local derivation if for every $x$, $y\in
A$, there exists a derivation $D_{x,y} : A\to A$ such that
$\bigtriangleup(x)=D_{x,y}(x)$ and $\bigtriangleup(y)=D_{x,y}(y)$.

In 1997, P. \v{S}emrl introduced the notion of 2-local derivations
and described 2-local derivations on the algebra $B(H)$ of all
bounded linear operators on the infinite-dimensional separable
Hilbert space $H$. A similar description for the
finite-dimensional case appeared later in (Kim and Kim 2004). In
the paper (Lin and Wong 2006) 2-local derivations have been
described on matrix algebras over finite-dimensional division
rings.

In (Ayupov and Kudaybergenov 2012) the authors suggested a new
technique and have generalized the above mentioned results of
(\v{S}emrl 1997) and (Kim and Kim 2004) for arbitrary Hilbert
spaces. Namely they considered 2-local derivations on the algebra
$B(H)$ of all linear bounded operators on an arbitrary (no
separability is assumed) Hilbert space $H$ and proved that every
2-local derivation on $B(H)$ is a derivation. After it is also
published a number of paper devoted to 2-local derivations on
associative algebras.

In the present paper we also suggest another technique and apply
to various associative algebras of infinite dimensional
matrix-valued functions on a compact. As a result we will have
that every 2-local derivation on such an algebra is a derivation.
As the main result of the paper it is established that every
2-local derivation on a $*$-algebra $C(Q, M_n(F))$ or
$C(Q,\mathcal{N}_n(F))$, where $Q$ is a compact, $M_n(F)$ is the
$*$-algebra of infinite dimensional matrices over complex numbers
(real numbers or quaternoins) (see section 1), $\mathcal{N}_n(F)$
is the $*$-algebra of infinite dimensional matrices over complex
numbers (real numbers or quaternoins) defined in section 2, is a
derivation. Also we explain that the method developed in the given
paper can be applied to Jordan and Lie algebras of infinite
dimensional matrix-valued functions on a compact.

We conclude that there are a number of various associative
algebras of infinite dimensional matrix-valued functions on a
compact every 2-local derivation of which is a derivation. The
main results of this paper are new and never proven. The method of
proving of these results represented in this paper is sufficiently
universal and can be applied to associative, Lie and Jordan
algebras. Its respective modification allows to prove similar
problem for Jordan and Lie algebras of infinite dimensional
matrix-valued functions on a compact.

\section{Preliminaries}

Let $M$ be an associative algebra.

{\it Definition.} A linear map $D : M\to M$ is called a
derivation, if $D(xy)=D(x)y+xD(y)$ for every two elements $x$,
$y\in M$.

A map $\Delta : M\to M$ is called a 2-local derivation, if for
every two elements $x$, $y\in M$ there exists a derivation
$D_{x,y}:M\to M$ such that $\Delta (x)=D_{x,y}(x)$, $\Delta
(y)=D_{x,y}(y)$.

It is known that each derivation $D$ on a von Neumann algebra $M$
is an inner derivation, that is there exists an element $a\in M$
such that
$$
D(x)=ax-xa, x\in M.
$$
Therefore for a von Neumann algebra $M$ the above definition is
equivalent to the following one: A map $\Delta : M\to M$ is called
a 2-local derivation, if for every two elements $x$, $y\in M$
there exists an element $a\in M$ such that $\Delta (x)=ax-xa$,
$\Delta (y)=ay-ya$.

Let throughout the paper $n$ be an arbitrary infinite cardinal
number, $\Xi$ be a set of indexes of the cardinality $n$. Let
$\{e_{ij}\}$ be a set of matrix units such that $e_{ij}$ is a
$n\times n$-dimensional matrix, i.e.
$e_{ij}=(a_{\alpha\beta})_{\alpha\beta\in\Xi}$, the $(i,j)$-th
component of which is $1$, i.e. $a_{ij}=1$, and the rest
components are zeros. Let $\{m_\xi\}_{\xi\in \Xi}$  be a set of
$n\times n$-dimensional matrixes. By $\sum_{\xi\in \Xi} m_\xi$ we
denote the matrix whose components are sums of the corresponding
components of matrixes of the set $\{m_\xi \}_{\xi\in \Xi}$. Let
throughout the paper $F={\mathbb C}$ (complexes), ${\mathbb
R}$(reals) or ${\mathbb H}$(quaternions) and
$$
M_n(F)=\{\{\lambda^{ij}e_{ij}\}: \,for\,\, all\,\, indices\,\,
i,\,j \,\lambda^{ij}\in F,
$$
$$
and\,\, there\,\, exists\,\, such\,\, number\,\, K\in {\bf
R},\,\,that \,\, for \,\, all\,\, n\in N
$$
$$
and\,\, \{e_{kl}\}_{kl=1}^n\subseteq \{e_{ij}\} \Vert\sum_{kl=1}^n
\lambda^{kl}e_{kl}\Vert \le K\},
$$
where $\Vert \,\, \Vert$ is a norm of a matrix. It is easy to see
that $M_n(F)$ is a vector space.

In the vector space
$$
\mathcal{M}_n(F)=\{\{\lambda^{ij}e_{ij}\}: \,for\,\, all\,\,
indices\,\, i,\,j \,\lambda^{ij}\in F\}
$$
of all $n\times n$-dimensional {\it matrices} (indexed sets) over
$F$ we introduce an associative multiplication as follows: if
$x=\{\lambda^{ij}e_{ij}\}$, $ y=\{\mu^{ij}e_{ij}\}$ are elements
of $\mathcal{M}_n(F)$ then $xy=\{ \sum_{\xi\in \Xi}
\lambda^{i\xi}\mu^{\xi j}e_{ij}\}$. With respect to this operation
$M_n(F)$ becomes an associative algebra and $M_n(F)\cong
B(l_2(\Xi))$, where $l_2(\Xi)$ is a Hilbert space over $F$ with
elements $\{x_i\}_{i\in \Xi}$, $x_i\in F$ for all $i\in \Xi$,
$B(l_2(\Xi))$ is the associative algebra of all bounded linear
operators on the Hilbert space $l_2(\Xi)$. Then $M_n(F)$ is a von
Neumann algebra of infinite $n\times n$-dimensional matrices over
$F$.

Similarly, if we take the algebra $B(H)$ of all bounded linear
operators on an arbitrary Hilbert space $H$ and if $\{q_i\}$ is an
arbitrary maximal orthogonal set of minimal projections of $B(H)$,
then $B(H)=\sum_{ij}^\oplus q_i B(H)q_j$ (see (Arzikulov 2008)).

Let throughout the paper $X$ be a hyperstonean compact, $C(X)$
denote the algebra of all $F$-valued continuous functions on $X$
and
$$
\mathcal{M}=\{\{\lambda^{ij}(x)e_{ij}\}_{ij\in\Xi}: (\forall
ij\,\,\, \lambda^{ij}(x)\in C(X))
$$
$$
(\exists K\in {\mathbb R})(\forall m\in N)(\forall
\{e_{kl}\}_{kl=1}^m\subseteq \{e_{ij}\})\Vert\sum_{kl=1\dots
m}\lambda^{kl}(x)e_{kl}\Vert\leq K\},
$$
where $\Vert\sum_{kl=1\dots m}\lambda^{kl}(x)e_{kl}\Vert\leq K$
means $(\forall x_o\in X) \Vert \sum_{kl=1\dots
m}\lambda^{kl}(x_o)e_{kl}\Vert\leq K$. The set $\mathcal{M}$ is a
vector space with point-wise algebraic operations. The map
$\Vert\,\,\, \Vert : \mathcal{M}\to {\bf R}_+$ defined as
$$
\Vert a \Vert = \sup_{\{e_{kl}\}_{kl=1}^n\subseteq
\{e_{ij}\}}\Vert\sum_{kl=1}^n \lambda^{kl}(x)e_{kl}\Vert,
$$
is a norm on the vector space $\mathcal{M}$, where $a\in
\mathcal{M}$ and $a=\{\lambda^{ij}(x)e_{ij}\}$.

In the vector space $V=\{\{\lambda^{ij}e_{ij}\}_{ij}:
\{\lambda^{ij}\}\subset C(X)\}$ of all infinite $n\times
n$-dimensional {\it matrices} (indexed sets) over $C(X)$ we
introduce an associative multiplication as follows: if
$x=\{\lambda^{ij}(x)e_{ij}\}$, $ y=\{\mu^{ij}(x)e_{ij}\}$ are
elements of $V$ then $xy=\{\sum_\xi\lambda^{i\xi}(x)\mu^{\xi
j}(x)e_{ij}\}$. With respect to this multiplication $\mathcal{M}$
becomes an associative algebra and $\mathcal{M}\cong C(X)\otimes
M_n(F)$. Thus $\mathcal{M}$ is a real or complex von Neumann
algebra of type I$_n$.

Let $M$ be a C$^*$-algebra, $\bigtriangleup :M\to M$ be a 2-local
derivation. Now let us show that $\bigtriangleup$ is homogeneous.
Indeed, for each $x\in M$, and for $\lambda \in {\mathbb C}$ there
exists a derivation $D_{x,\lambda x}$ such that
$\bigtriangleup(x)=D_{x,\lambda x}(x)$ and $\bigtriangleup(\lambda
x)=D_{x,\lambda x}(\lambda x)$. Then
$$
\bigtriangleup(\lambda x)=D_{x,\lambda x}(\lambda x)=\lambda D_{x,\lambda x}(x) =\lambda \bigtriangleup(x).
$$
Hence, $\bigtriangleup$ is homogenous. At the same time, for each
$x\in M$, there exists a derivation   $D_{x,x^2}$ such that
$\bigtriangleup(x)=D_{x,x^2}(x)$ and
$\bigtriangleup(x^2)=D_{x,x^2}(x^2)$. Then
$$
\bigtriangleup(x^2)=D_{x,x^2}(x^2)=D_{x,x^2}(x)x+xD_{x,x^2}(x) =\bigtriangleup(x)x+x\bigtriangleup(x).
$$

In (Bresar 1988) it is proved that every Jordan derivation on a
semi-prime algebra is a derivation. Since $M$ is semi-prime (i.e.
$aMa = \{0\}$ implies that $a = \{0\}$), the map $\bigtriangleup$
is a derivation if it is additive. Therefore, to prove that the
2-local derivation $\bigtriangleup :M\to M$ is a derivation it is
sufficient to prove that $\bigtriangleup :M\to M$ is additive in
the proof of theorem 1.

\section{2-local derivations on some associative algebras of matrix-valued
functions}

Let $Q$ be a compact. Then there exists a hyperstonean compact $X$
such that for the algebra $C(Q)$ of all continuous complex
number-valued functions on $Q$ we have $C(Q)^{**}\cong C(X)$. If
we take the $*$-algebra $C(Q, M_n({\mathbb C}))$ of all continuous
maps of $Q$ to $M_n({\mathbb C})$, then we may assume that $C(Q,
M_n({\mathbb C}))\subseteq \mathcal{M}$. In this case the set
$\{e_{ij}\}$ of constant functions belongs to $C(Q, M_n({\mathbb
C}))$ and the weak closure of $C(Q, M_n({\mathbb C}))$ in
$\mathcal{M}$ coincides with $\mathcal{M}$. Hence by separately
weakly continuity of multiplication every derivation of $C(Q,
M_n({\mathbb C}))$ has a unique extension to a derivation on
$\mathcal{M}$. Therefore, if $\bigtriangleup$ is a 2-local
derivation on $C(Q, M_n({\mathbb C}))$, then for every two
elements $x$, $y\in C(Q, M_n({\mathbb C}))$ there exists a
derivation $D_{x,y}:\mathcal{M}\to \mathcal{M}$ such that
$\bigtriangleup (x)=D_{x,y}(x)$, $\bigtriangleup (y)=D_{x,y}(y)$,
i.e. $D_{x,y}$ is a derivation of $\mathcal{M}$ (not only of $C(Q,
M_n({\mathbb C}))$). The following theorem is the key  result of
this section.

{\bf Theorem 1.} {\it Let $\bigtriangleup$ be a 2-local derivation
on $C(Q, M_n({\mathbb C}))$. Then $\bigtriangleup$ is a
derivation.}

\medskip

First let us prove lemmata which are necessary for the proof of
theorem 1.

By the above arguments for every 2-local derivation
$\bigtriangleup$ on $C(Q, M_n(F))$ and for each $x\in C(Q,
M_n(F))$ there exist $a\in \mathcal{M}$ such that
$$
\bigtriangleup(x)=ax-xa.
$$
Put
$$
e_{ij}:=\{\lambda^{\xi\eta}e_{\xi\eta}\},
$$
where for all $\xi$, $\eta$, if $\xi=i$, $\eta=j$ then
$\lambda^{\xi\eta}={\bf 1}$, else $\lambda^{\xi\eta}=0$, ${\bf 1}$
is unit of the algebra $C(Q)$. Let $\{a(ij)\}\subset \mathcal{M}$
be a subset such that
$$
\bigtriangleup(e_{ij})=a(ij)e_{ij}-e_{ij}a(ij).
$$
for all $i$, $j$,  let $a^{ij}e_{ij}$, $a^{ij}\in C(Q)$, be the
(i,j)-th component of the element $e_{ii}a(ji)e_{jj}$ of
$\mathcal{M}$ for all pairs of different indexes $i$, $j$ and let
$\{a^{\xi\eta}e_{\xi\eta}\}_{\xi\neq \eta}$ be the matrix in $V$
with all such components, the diagonal components of which are
zeros.

\medskip

{\bf Lemma 2.} {\it For each pair $i$, $j$ of different indices
the following equality is valid
$$
\bigtriangleup(e_{ij})=\{a^{\xi\eta}e_{\xi\eta}\}_{\xi\neq
\eta}e_{ij}-e_{ij}\{a^{\xi\eta}e_{\xi\eta}\}_{\xi\neq \eta}+
a(ij)_{ii}e_{ij}-e_{ij}a(ij)_{jj},    \,\,\,\,\,\,\,(1)
$$
where $a(ij)^{ii}$, $a(ij)^{jj}$ are functions in $C(Q)$ which are
the coefficients of the Peirce components $e_{ii}a(ij)e_{ii}$,
$e_{jj}a(ij)e_{jj}$.}

{\it Proof.} Let $k$ be an arbitrary index different from $i$, $j$
and let $a(ij,ik)\in \mathcal{M}$ be an element such that
$$
\bigtriangleup(e_{ik})=a(ij,ik)e_{ik}-e_{ik}a(ij,ik) \,\, \text{and}\,\,
\bigtriangleup(e_{ij})=a(ij,ik)e_{ij}-e_{ij}a(ij,ik).
$$
Then
$$
e_{kk}\bigtriangleup(e_{ij})e_{jj}=e_{kk}(a(ij,ik)e_{ij}-e_{ij}a(ij,ik))e_{jj}=
$$
$$
e_{kk}a(ij,ik)e_{ij}-0=e_{kk}a(ik)e_{ij}-e_{kk}e_{ij}\{a^{\xi\eta}e_{\xi\eta}\}_{\xi\neq
\eta}e_{jj}=
$$
$$
e_{kk}a_{ki}e_{ij}-e_{kk}e_{ij}\{a^{\xi\eta}e_{\xi\eta}\}_{\xi\neq
\eta}e_{jj}= e_{kk}\{a^{\xi\eta}e_{\xi\eta}\}_{\xi\neq
\eta}e_{ij}-e_{kk}e_{ij}\{a^{\xi\eta}e_{\xi\eta}\}_{\xi\neq
\eta}e_{jj}=
$$
$$
e_{kk}(\{a^{\xi\eta}e_{\xi\eta}\}_{\xi\neq
\eta}e_{ij}-e_{ij}\{a^{\xi\eta}e_{\xi\eta}\}_{\xi\neq
\eta})e_{jj}.
$$

Similarly,
$$
e_{kk}\bigtriangleup(e_{ij})e_{ii}=
e_{kk}(\{a^{\xi\eta}e_{\xi\eta}\}_{\xi\neq
\eta}e_{ij}-e_{ij}\{a^{\xi\eta}e_{\xi\eta}\}_{\xi\neq
\eta})e_{ii}.
$$

Let $a(ij,kj)\in \mathcal{M}$ be an element such that
$$
\bigtriangleup(e_{kj})=a(ij,kj)e_{kj}-e_{kj}a(ij,kj)  \,\, \text{and}\,\,
\bigtriangleup(e_{ij})=a(ij,kj)e_{ij}-e_{ij}a(ij,kj).
$$

Then
$$
e_{ii}\bigtriangleup(e_{ij})e_{kk}=e_{ii}(a(ij,kj)e_{ij}-e_{ij}a(ij,kj))e_{kk}=
$$
$$
0-e_{ij}a(ij,kj)e_{kk}=0-e_{ij}a(kj)e_{kk}=0-e_{ij}a_{jk}e_{kk}=
$$
$$
e_{ii}\{a^{\xi\eta}e_{\xi\eta}\}_{\xi\neq
\eta}e_{ij}e_{kk}-e_{ij}\{a^{\xi\eta}e_{\xi\eta}\}_{\xi\neq
\eta}e_{kk}=
$$
$$
e_{ii}(\{a^{\xi\eta}e_{\xi\eta}\}_{\xi\neq
\eta}e_{ij}-e_{ij}\{a^{\xi\eta}e_{\xi\eta}\}_{\xi\neq
\eta})e_{kk}.
$$

Also similarly we have
$$
e_{jj}\bigtriangleup(e_{ij})e_{kk}=
e_{jj}(\{a^{\xi\eta}e_{\xi\eta}\}_{\xi\neq
\eta}e_{ij}-e_{ij}\{a^{\xi\eta}e_{\xi\eta}\}_{\xi\neq
\eta})e_{kk},
$$
$$
e_{ii}\bigtriangleup(e_{ij})e_{ii}=
e_{ii}(\{a^{\xi\eta}e_{\xi\eta}\}_{\xi\neq
\eta}e_{ij}-e_{ij}\{a^{\xi\eta}e_{\xi\eta}\}_{\xi\neq
\eta})e_{ii},
$$
$$
e_{jj}\bigtriangleup(e_{ij})e_{jj}=
e_{jj}(\{a^{\xi\eta}e_{\xi\eta}\}_{\xi\neq
\eta}e_{ij}-e_{ij}\{a^{\xi\eta}e_{\xi\eta}\}_{\xi\neq
\eta})e_{jj}.
$$

Hence the equality (1) is valid. $\triangleright$

We take elements of the sets $\{\{e_{i\xi}\}_\xi\}_i$ and $\{\{e_{\xi j}\}_\xi\}_j$ in
pairs $(\{e_{\alpha\xi}\}_\xi,\{e_{\xi \beta}\}_\xi)$ such that $\alpha\neq \beta$.
Then using the set $\{(\{e_{\alpha\xi}\}_\xi,\{e_{\xi\beta}\}_\xi)\}$ of such pairs
we get the set $\{e_{\alpha\beta}\}$.

Let $x_o=\{e_{\alpha\beta}\}$ be a set $\{v^{ij}e_{ij}\}_{ij}$
such that for all $i$, $j$ if $(\alpha,\beta)\neq (i,j)$ then
$v_{ij}=0\in C(Q)$ else $v_{ij}={\bf 1}\in C(Q)$. Then $x_o\in
C(Q, M_n({\mathbb C}))$. Fix different indices $i_o$, $j_o$. Let
$c\in \mathcal{M}$ be an element such that
$$
\bigtriangleup(e_{i_oj_o})=ce_{i_oj_o}-e_{i_oj_o}c \,\,
\text{and}\,\, \bigtriangleup(x_o)=cx_o-x_oc.
$$

Put $c=\{c^{ij}e_{ij}\}\in \mathcal{M}$ and
$\bar{a}=\{a^{ij}e_{ij}\}_{i\neq j}\cup \{a^{ii}e_{ii}\}$, where
$\{a^{ii}e_{ii}\}=\{c^{ii}e_{ii}\}$.

\medskip

{\bf Lemma 3.} {\it Let $\xi$, $\eta$ be arbitrary different
indices, and let $b=\{b^{ij}e_{ij}\}\in \mathcal{M}$ be an element
such that
$$
\bigtriangleup(e_{\xi\eta})=be_{\xi\eta}-e_{\xi\eta}b \,\,
\text{and}\,\, \bigtriangleup(x_o)=bx_o-x_ob.
$$
Then $c^{\xi\xi}-c^{\eta\eta}=b^{\xi\xi}-b^{\eta\eta}$.}

{\it Proof.} We have that there exist $\bar{\alpha}$, $\bar{\beta}$ such that
$e_{\xi\bar{\alpha}}$, $e_{\bar{\beta}\eta}\in \{e_{\alpha\beta}\}$ (or $e_{\bar{\alpha}\eta}$,
$e_{\xi\bar{\beta}}\in \{e_{\alpha\beta}\}$, or $e_{\bar{\alpha},\bar{\beta}}\in \{e_{\alpha\beta}\}$), and there
exists a chain of pairs of indexes $(\hat{\alpha}, \hat{\beta})$ in $\Omega$, where
$\Omega=\{(\check{\alpha},\check{\beta}): e_{\check{\alpha},\check{\beta}}\in \{e_{\alpha\beta}\}\}$,
connecting pairs $(\xi, \bar{\alpha})$, $(\bar{\beta},\eta)$ i.e.,
$$
(\xi,\bar{\alpha}), (\bar{\alpha}, \xi_1), (\xi_1, \eta_1),\dots ,(\eta_2, \bar{\beta}), (\bar{\beta},\eta).
$$
Then
$$
c^{\xi\xi}-c^{\bar{\alpha}\bar{\alpha}}=b^{\xi\xi}-b^{\bar{\alpha}\bar{\alpha}},
c^{\bar{\alpha}\bar{\alpha}}-c^{\xi_1\xi_1}=b^{\bar{\alpha}\bar{\alpha}}-b^{\xi_1\xi_1},
$$
$$
c^{\xi_1\xi_1}-c^{\eta_1\eta_1}=b^{\xi_1\xi_1}-b^{\eta_1\eta_1},\dots
,
c^{\eta_2\eta_2}-c^{\bar{\beta}\bar{\beta}}=b^{\eta_2\eta_2}-b^{\bar{\beta}\bar{\beta}},
c^{\bar{\beta}\bar{\beta}}-c^{\eta\eta}=b^{\bar{\beta}\bar{\beta}}-b^{\eta\eta}.
$$
Hence
$$
c^{\xi\xi}-b^{\xi\xi}=c^{\bar{\alpha}\bar{\alpha}}-b^{\bar{\alpha}\bar{\alpha}},
c^{\bar{\alpha}\bar{\alpha}}-b^{\bar{\alpha}\bar{\alpha}}=c^{\xi_1\xi_1}-b^{\xi_1\xi_1},
$$
$$
c^{\xi_1\xi_1}-b^{\xi_1\xi_1}=c^{\eta_1\eta_1}-b^{\eta_1\eta_1},\dots
,
c^{\eta_2\eta_2}-b^{\eta_2\eta_2}=c^{\bar{\beta}\bar{\beta}}-b^{\bar{\beta}\bar{\beta}},
c^{\bar{\beta}\bar{\beta}}-b^{\bar{\beta}\bar{\beta}}=c^{\eta\eta}-b^{\eta\eta}.
$$
and $c^{\xi\xi}-b^{\xi\xi}=c^{\eta\eta}-b^{\eta\eta}$,
$c^{\xi\xi}-c^{\eta\eta}=b^{\xi\xi}-b^{\eta\eta}$.

Therefore $c^{\xi\xi}-c^{\eta\eta}=b^{\xi\xi}-b^{\eta\eta}$.
$\triangleright$

\medskip

{\bf Lemma 4.} {\it Let $x$ be an element of the algebra $C(Q,
M_n({\mathbb C}))$. Then
$$
\bigtriangleup(x)=\bar{a}x-x\bar{a},
$$
where $\bar{a}$ is defined as above.}

{\it Proof.}
Let $d(ij)\in \mathcal{M}$ be an element such that
$$
\bigtriangleup(e_{ij})=d(ij)e_{ij}-e_{ij}d(ij) \,\, \text{and}\,\,
\bigtriangleup(x)=d(ij)x-xd(ij)
$$
and $i\neq j$. Then
$$
\bigtriangleup(e_{ij})=d(ij)e_{ij}-e_{ij}d(ij)=
$$
$$
e_{ii}d(ij)e_{ij}-e_{ij}d(ij)e_{jj}+
(1-e_{ii})d(ij)e_{ij}-e_{ij}d(ij)(1-e_{jj})=
$$
$$
a(ij)_{ii}e_{ij}-e_{ij}a(ij)_{jj}+
\{a^{\xi\eta}e_{\xi\eta}\}_{\xi\neq
\eta}e_{ij}-e_{ij}\{a^{\xi\eta}e_{\xi\eta}\}_{\xi\neq \eta}
$$
for all $i$, $j$ by lemma 2.

Since
$$
e_{ii}d(ij)e_{ij}-e_{ij}d(ij)e_{jj}=a(ij)_{ii}e_{ij}-e_{ij}a(ij)_{jj}
$$
we have
$$
(1-e_{ii})d(ij)e_{ii}=\{a^{\xi\eta}e_{\xi\eta}\}_{\xi\neq
\eta}e_{ii},
$$
$$
e_{jj}d(ij)(1-e_{jj})=e_{jj}\{a^{\xi\eta}e_{\xi\eta}\}_{\xi\neq
\eta}
$$
for all different $i$ and $j$.

Let $b=\{b^{ij}e_{ij}\}\in \mathcal{M}$ be an element such that
$$
\bigtriangleup(e_{ij})=be_{ij}-e_{ij}b \,\, \text{and}\,\,
\bigtriangleup(x_o)=bx_o-x_ob.
$$
Then $b^{ii}-b^{jj}=c^{ii}-c^{jj}$ by lemma 3. We have
$b^{ii}-b^{jj}=d(ij)^{ii}-d(ij)^{jj}$ since
$$
be_{ij}-e_{ij}b=d(ij)e_{ij}-e_{ij}d(ij).
$$
Hence
$$
c^{ii}-c^{jj}=d(ij)^{ii}-d(ij)^{jj},
c^{jj}-c^{ii}=d(ij)^{jj}-d(ij)^{ii}.
$$

Therefore we have
$$
e_{jj}\bigtriangleup(x)e_{ii}=e_{jj}(d(ij)x-xd(ij))e_{ii}=
$$
$$
e_{jj}d(ij)(1-e_{jj})xe_{ii}+
e_{jj}d(ij)e_{jj}xe_{ii}-e_{jj}x(1-e_{ii})d(ij)e_{ii}-e_{jj}xe_{ii}d(ij)e_{ii}=
$$
$$
e_{jj}\{a^{\xi\eta}e_{\xi\eta}\}_{\xi\neq
\eta}xe_{ii}-e_{jj}x\{a^{\xi\eta}e_{\xi\eta}\}_{\xi\neq
\eta}e_{ii}+ e_{jj}d(ij)e_{jj}xe_{ii}-e_{jj}xe_{ii}d(ij)e_{ii}=
$$
$$
e_{jj}\{a^{\xi\eta}e_{\xi\eta}\}_{\xi\neq
\eta}xe_{ii}-e_{jj}x\{a^{\xi\eta}e_{\xi\eta}\}_{\xi\neq
\eta}e_{ii}+ c^{jj}e_{jj}xe_{ii}-e_{jj}xe_{ii}c^{ii}e_{ii}=
$$
$$
e_{jj}\{a^{\xi\eta}e_{\xi\eta}\}_{\xi\neq
\eta}xe_{ii}-e_{jj}x\{a^{\xi\eta}e_{\xi\eta}\}_{\xi\neq
\eta}e_{ii}+
$$
$$
e_{jj}(\sum_\xi a^{\xi\xi}e_{\xi\xi})xe_{ii}-e_{jj}x(\sum_\xi
a^{\xi\xi}e_{\xi\xi})e_{ii}=
$$
$$
e_{jj}\{a^{\xi\eta}e_{\xi\eta}\}xe_{ii}-e_{jj}x\{a^{\xi\eta}e_{\xi\eta}\}e_{ii}=
e_{jj}(\bar{a}x-x\bar{a})e_{ii}.
$$

Let $d(ii)$, $v$, $w\in \mathcal{M}$ be elements such that
$$
\bigtriangleup(e_{ii})=d(ii)e_{ii}-e_{ii}d(ii) \,\, \text{and}\,\,
\bigtriangleup(x)=d(ii)x-xd(ii),
$$
$$
\bigtriangleup(e_{ii})=ve_{ii}-e_{ii}v,
\bigtriangleup(e_{ij})=ve_{ij}-e_{ij}v,
$$
and
$$
\bigtriangleup(e_{ii})=we_{ii}-e_{ii}w,
\bigtriangleup(e_{ji})=we_{ji}-e_{ji}w.
$$
Then
$$
(1-e_{ii})a(ij)e_{ii}=(1-e_{ii})ve_{ii}=(1-e_{ii})d(ii)e_{ii},
$$
and
$$
e_{ii}a(ji)(1-e_{ii})=e_{ii}w(1-e_{ii})=e_{ii}d(ii)(1-e_{ii}).
$$
By lemma 2
$$
\bigtriangleup(e_{ij})=a(ij)e_{ij}-e_{ij}a(ij)=
$$
$$
\{a^{\xi\eta}e_{\xi\eta}\}_{\xi\neq
\eta}e_{ij}-e_{ij}\{a^{\xi\eta}e_{\xi\eta}\}_{\xi\neq \eta}+
a(ij)^{ii}e_{ij}-e_{ij}a(ij)^{jj}
$$
and
$$
(1-e_{ii})a(ij)e_{ii}=\{a^{\xi\eta}e_{\xi\eta}\}_{\xi\neq
\eta}e_{ii}.
$$
Similarly
$$
e_{ii}a(ji)(1-e_{ii})=e_{ii}\{a^{\xi\eta}e_{\xi\eta}\}_{\xi\neq
\eta}.
$$

Hence
$$
e_{ii}\bigtriangleup(x)e_{ii}=e_{ii}(d(ii)x-xd(ii))e_{ii}=
$$
$$
e_{ii}d(ii)(1-e_{ii})xe_{ii}+
e_{ii}d(ii)e_{ii}xe_{ii}-e_{ii}x(1-e_{ii})d(ii)e_{ii}-e_{ii}xe_{ii}d(ii)e_{ii}=
$$
$$
e_{ii}a(ji)(1-e_{ii})xe_{ii}+
e_{ii}d(ii)e_{ii}xe_{ii}-e_{ii}x(1-e_{ii})a(ij)e_{ii}-e_{ii}xe_{ii}d(ii)e_{ii}=
$$
$$
e_{ii}\{a^{\xi\eta}e_{\xi\eta}\}_{\xi\neq
\eta}xe_{ii}-e_{ii}x\{a^{\xi\eta}e_{\xi\eta}\}_{\xi\neq
\eta}e_{ii}+ e_{ii}d(ii)e_{ii}xe_{ii}-e_{ii}xe_{ii}d(ii)e_{ii}=
$$
$$
e_{ii}\{a^{\xi\eta}e_{\xi\eta}\}_{\xi\neq
\eta}xe_{ii}-e_{ii}x\{a^{\xi\eta}e_{\xi\eta}\}_{\xi\neq
\eta}e_{ii}+c^{ii}e_{ii}xe_{ii}-e_{ii}xc^{ii}e_{ii}=
$$
$$
e_{ii}\{a^{\xi\eta}e_{\xi\eta}\}_{\xi\neq
\eta}xe_{ii}-e_{ii}x\{a^{\xi\eta}e_{\xi\eta}\}_{\xi\neq
\eta}e_{ii}+
$$
$$
e_{ii}(\sum_\xi a^{\xi\xi}e_{\xi\xi})xe_{ii}-e_{ii}x(\sum_\xi
a^{\xi\xi}e_{\xi\xi})e_{ii}=
$$
$$
e_{ii}\{a^{\xi\eta}e_{\xi\eta}\}xe_{ii}-e_{ii}x\{a^{\xi\eta}e_{\xi\eta}\}e_{ii}=e_{ii}(\bar{a}x-x\bar{a})e_{ii}.
$$

Hence
$$
\bigtriangleup(x)=\bar{a}x-x\bar{a}
$$
for all $x\in C(Q, M_n({\mathbb C}))$. $\triangleright$

\medskip

{\it Proof of theorem 1.} By lemma 4
$\bigtriangleup(e_{ii})=\bar{a}e_{ii}-e_{ii}\bar{a}\in
\mathcal{M}$. Hence
$$
\sum_\xi a^{\xi i}e_{\xi i}-\sum_\xi a^{i\xi}e_{i\xi}\in
\mathcal{M}.
$$
Then
$$
e_{ii}(\sum_\xi a^{\xi i}e_{\xi i}-\sum_\xi a^{i\xi}e_{i\xi})=
a^{ii}e_{ii}-\sum_\xi a^{i\xi}e_{i\xi}\in \mathcal{M}
$$
and
$$
(\sum_\xi a^{\xi i}e_{\xi i}-\sum_\xi a^{i\xi}e_{i\xi})e_{ii}=
\sum_\xi a^{\xi i}e_{\xi i}-a^{ii}e_{ii}\in \mathcal{M}.
$$
Therefore $\sum_\xi a^{\xi i}e_{\xi i}$, $\sum_\xi
a^{i\xi}e_{i\xi}\in \mathcal{M}$ i.e., $\bar{a}e_{ii},
e_{ii}\bar{a}\in \mathcal{M}$. Hence $e_{ii}\bar{a}x,
x\bar{a}e_{ii}\in \mathcal{M}$ for each $i$ and
$$
\bar{a}x, x\bar{a}\in V
$$
for each element $x=\{x^{ij}e_{ij}\}\in C(Q, M_n({\mathbb C}))$,
i.e.,
$$
\sum_\xi a^{i\xi}x^{\xi j}e_{ij}, \sum_\xi x^{i\xi}a^{\xi
j}e_{ij}\in C(Q)e_{ij}
$$
for all $i$, $j$. Therefore for all $x$, $y\in C(Q, M_n({\mathbb
C}))$ we have that the elements $\bar{a}x$, $x\bar{a}$,
$\bar{a}y$, $y\bar{a}$, $\bar{a}(x+y)$, $(x+y)\bar{a}$ belong to
$V$. Hence
$$
\bigtriangleup(x+y)=\bigtriangleup(x)+\bigtriangleup(y)
$$
by lemma 4.

Similarly for all $x$, $y\in C(Q, M_n({\mathbb C}))$ we have
$$
(\bar{a}x+x\bar{a})y=\bar{a}xy-x\bar{a}y\in \mathcal{M}, \bar{a}xy=\bar{a}(xy)\in V.
$$
Then $x\bar{a}y=\bar{a}xy-(\bar{a}x-x\bar{a})y$ and $x\bar{a}y\in V$.
Therefore
$$
\bar{a}(xy)-(xy)\bar{a}=\bar{a}xy-x\bar{a}y+x\bar{a}y-xy\bar{a}=(\bar{a}x-x\bar{a})y+x(\bar{a}y-y\bar{a}).
$$
Hence
$$
\bigtriangleup(xy)=\bigtriangleup(x)y+x\bigtriangleup(y)
$$
by lemma 4. By section 1 $\bigtriangleup$ is homogeneous. Hence,
$\bigtriangleup$ is a linear operator and a derivation. The proof
is complete. $\triangleright$

If we take the $*$-algebra $C(Q, M_n(F))$, $F={\mathbb R}$ or
${\mathbb H}$, then we can similarly prove the following theorem.

{\bf Theorem 2.} {\it Let $\bigtriangleup$ be a 2-local derivation
on $C(Q, M_n(F))$. Then $\bigtriangleup$ is a derivation.}

Note, for theorem 2 to be proved the proof of theorem 1 will be
repeated with very minor modification.

Let
$$
\sum_{ij}^o Fe_{ij}=\{\{\lambda^{ij}e_{ij}\}: \,for\,\, all\,\,
indices\,\, i,\,j \,\lambda^{ij}\in F,\,and\,
$$
$$
\forall \varepsilon>0 \exists n_o\in N \,\,such\,\,that\,\,\forall
n\geq m\geq n_o
$$
$$
\Vert
\sum_{i=m}^n[\sum_{k=1,\dots,i-1}(\lambda^{ki}e_{ki}+\lambda^{ik}e_{ik})+\lambda^{ii}e_{ii}]\Vert<\varepsilon\}.
$$
where $\Vert \,\, \Vert$ is a norm of a matrix. Then $\sum_{ij}^o
Fe_{ij}$ is a C$^*$-algebra with respect to componentwise
algebraic operations, the bilinear operation and the norm
(Arzikulov 2012). Since $Fe_{ij}$ is a simple C$^*$-algebra for
all $i$, then by the proof of theorem 8 in (Arzikulov 2012) the
C$^*$-algebra $\sum_{ij}^o Fe_{ij}$ is simple. Let
$\mathcal{N}_n(F)=\sum_{ij}^o Fe_{ij}$. Then
$C(Q,\mathcal{N}_n(F))$ is a real or complex C$^*$-algebra, where
$(F={\mathbb C}$, ${\mathbb R}$ or ${\mathbb H}$) and $C(Q,
\mathcal{N}_n(F))\subseteq \mathcal{M}$. Hence similar to theorems
1, 2 we can prove the following theorem

{\bf Theorem 3.} {\it Let $\bigtriangleup$ be a 2-local derivation
on $C(Q,\mathcal{N}_n(F))$. Then $\bigtriangleup$ is a
derivation.}

It is known that the set $\mathcal{M}_{sa}$ of all self-adjoint
elements (i.e. $a^*=a$) of $\mathcal{M}$ forms a Jordan algebra
with respect to the operation of multiplication $a\cdot
b=\frac{1}{2}(ab+ba)$. The following problem can be similarly
solved.

{\bf Problem 1.} Develop a Jordan analog of the method applied in
the proof of theorem 1 and prove that every 2-local derivation
$\bigtriangleup$ on the Jordan algebra $\mathcal{M}_{sa}$ or $C(Q,
M_n(F)_{sa})$ or $C(Q,\mathcal{N}_n(F)_{sa})$ is a derivation.

It is known that the set $\mathcal{M}_k=\{a\in \mathcal{M}:
a^*=-a\}$ forms a Lie algebra with respect to the operation of
multiplication $[a, b]=ab-ba$. So it is natural to consider  the following problem.

{\bf Problem 2.} Develop a Lie analog of the method applied in the
proof of theorem 1 and prove that every 2-local derivation
$\bigtriangleup$ on the Lie algebra $\mathcal{M}_k$ or $C(Q,
M_n(F)_k)$ or $C(Q,\mathcal{N}_n(F)_k)$ is a derivation.

\bigskip

{\bf Acknowledgements}

The authors want to thank K.K.Kudaybergenov for many stimulating
conversations on the subject.

\bigskip

{\bf References}

\bigskip

\v{S}emrl P (1997) Local automorphisms and derivations on $B(H)$.
Proc Amer Math Soc 125:2677-2680

Kim S O, Kim J S (2004) Local automorphisms and derivations on
$M_n$. Proc Amer Math Soc 132:1389-1392

Lin Y, Wong T (2006) A note on 2-local maps. Proc Edinb Math Soc
49:701-708.

Ayupov Sh A, Kudaybergenov K K (2012) 2-local derivations and
automorphisms on $B(H)$ 395:15-18

Arzikulov F N (2008) Infinite order and norm decompositions of
C*-algebras. Int Journal of Math Analysis  2(5):255-262

Arzikulov F N (2012) Infinite norm decompositions of C*-algebras.
In: Operator Theory: Advances and Applications, vol 220. Springer
Basel AG pp. 11-21

Bresar M (1988) Jordan derivations on semiprime rings. Proc Amer
Math Soc 104:1003-1006

\end{document}